\documentclass[10pt]{amsart}

\usepackage[utf8]{inputenc}
\usepackage{amssymb,tikz}
\usepackage{amsthm}
\usepackage{amsmath}
\usepackage{amsbsy, a4wide}
\usepackage[all]{xy}
\usepackage{bm}
\usepackage{hyperref,caption}

\begin{document}

\title[]{A hidden Signal in the Ulam sequence }
\author[]{Stefan Steinerberger}
\address{Department of Mathematics, Yale University, 10 Hillhouse Avenue, New Haven, CT 06511, USA}
\email{stefan.steinerberger@yale.edu}
\begin{abstract} The Ulam sequence is defined as $a_1 =1, a_2 = 2$ and $a_n$ being the smallest integer that can be written as the sum of
two distinct earlier elements in a unique way. This gives $$1, 2, 3, 4, 6, 8, 11, 13, 16, 18, 26, 28, 36, 38, 47, \dots$$ 
Ulam remarked that understanding the sequence, which has been described as 'quite erratic', seems difficult and indeed nothing 
is known. We report the empirical discovery of a surprising global rigidity phenomenon: there seems to exist a real $\alpha \sim 2.5714474995\dots$ such that 
$$\left\{\alpha a_n ~\mbox{mod}~2\pi : n\in \mathbb{N}\right\} \quad \mbox{generates an absolutely continuous \textit{non-uniform} measure}$$
supported on a subset of $\mathbb{T}$. Indeed, for the first $10^7$ elements of Ulam's sequence, 
$$ \cos{\left( 2.5714474995~ a_n\right)} < 0 \qquad \mbox{for all}~a_n \notin \left\{2, 3, 47, 69\right\}.$$
The same phenomenon arises for some other initial conditions $a_1, a_2$: the
distribution functions look very different from each other and have curious shapes. A similar but
more subtle phenomenon seems to arise in Lagarias' variant of MacMahon's 'primes of measurement' sequence.
\end{abstract}
\maketitle
\vspace{-25pt}

\section{Introduction}
Stanis\l{}aw Ulam introduced his sequence
 $$1, 2, 3, 4, 6, 8, 11, 13, 16, 18, 26, 28, 36, 38, 47, 48, 53, 57, 62, 69, 72, 77, 82, 87, 97 \dots \dots$$
in a 1964 survey \cite{ulam} on unsolved problems. The construction is given by $a_1 =1, a_2 = 2$ and then iteratively
choosing the smallest integer that can be written as the sum of two distinct earlier elements in a unique way as the next element. Ulam writes
\begin{quote}
One can consider a rule for growth of patterns -- in one dimension it would be
merely a rule for obtaining successive integers. [...] In both cases simple questions that come to
mind about the properties of a sequence of integers thus obtained are notoriously
hard to answer. (Ulam, 1964)
\end{quote}
He asks in the Preface of \cite{ulam2} whether it is possible to determine the asymptotic density of the sequence (this problem is
sometimes incorrectly attributed to Recaman \cite{recaman}). The Ulam sequence was first extensively computed in 1966 by Muller \cite{mull, wunderlich}, who gave the first 20.000 terms and predicted density 0. The most extensive computation we had access to is due to Daniel Strottman \cite{strott1}, who computed the first 11.172.261 elements (up to 150.999.995). This data set shows that some of Muller's predictions are not correct (the density, for example, seems to be very stable and around 0.0739). Different initial values $a_1, a_2$ can give rise to more structured sequences \cite{finch1, finch2, q}:
for some of them the sequence of consecutive differences $a_{n+1} - a_{n}$ is eventually periodic. It seems that this is not the case for Ulam's
sequence: Knuth \cite{oeis} remarks that $a_{4953}-a_{4952}=262$ and $a_{18858} - a_{18857} = 315$. The sequence 'does not appear to follow any recognizable pattern' \cite{finch3} and is 'quite erratic' \cite{schm}.
We describe the (accidental) discovery of some very surprising structure. While using Fourier series with Ulam numbers as frequencies, we
noticed a persisting signal in the noise: indeed, dilating the sequence by a factor $\alpha \sim 2.5714474995\dots$ and considering
the sequence $(\alpha a_n$ mod $2\pi)$ gives rise to a very regular distribution function. One surprising implication is
$$ \cos{\left(2.5714474995 ~a_n\right)} < 0 \qquad \mbox{for the first $10^7$ Ulam numbers except}~~\quad~2, 3, 47, 69.$$
The dilation factor $\alpha$ seems to be a universal constant (of which we were able to compute the first few digits); given the
distribution function, there is nothing special about the cosine and many similar inequalities could be derived (we used the cosine
because it is perhaps the simplest).
\section{The Observation}
\subsection{A Fourier series.} If the Ulam sequence had positive density, then by a classical theorem of
Roth \cite{roth} it would contain many 3-arithmetic progressions. We were interested in whether the existence of progressions would have any
impact on the structure of the sequence (because,  if $a\cdot n + b$ is in the sequence for $n=1,2,3$, then $a$ or $2a$ is not).  
There is a well-known connection between randomness in sets and smallness of Fourier coefficients and this motivated us to look at
$$ f_N(x) = \mbox{Re} \sum_{n=1}^{N}{e^{i a_n x} } = \sum_{n=1}^{N}{\cos{(a_n x)}}.$$
Clearly, $f_N(0) = N$ and $\|f_N\|_{L^2} \sim \sqrt{N}$ and therefore we expected $|f_N(x)| \sim \sqrt{N}$ for $x$ outside of 0. Much to our surprise, we discovered that this is not the case and that there is an $\alpha \sim 2.571\dots$ with
that $f_N(\alpha) \sim -0.8 N$. Such a behavior is, of course, an indicator of an embedded signal.

\begin{figure}[h!]
\begin{minipage}[t]{0.49\columnwidth}%
\begin{center}
\includegraphics[width = 5cm]{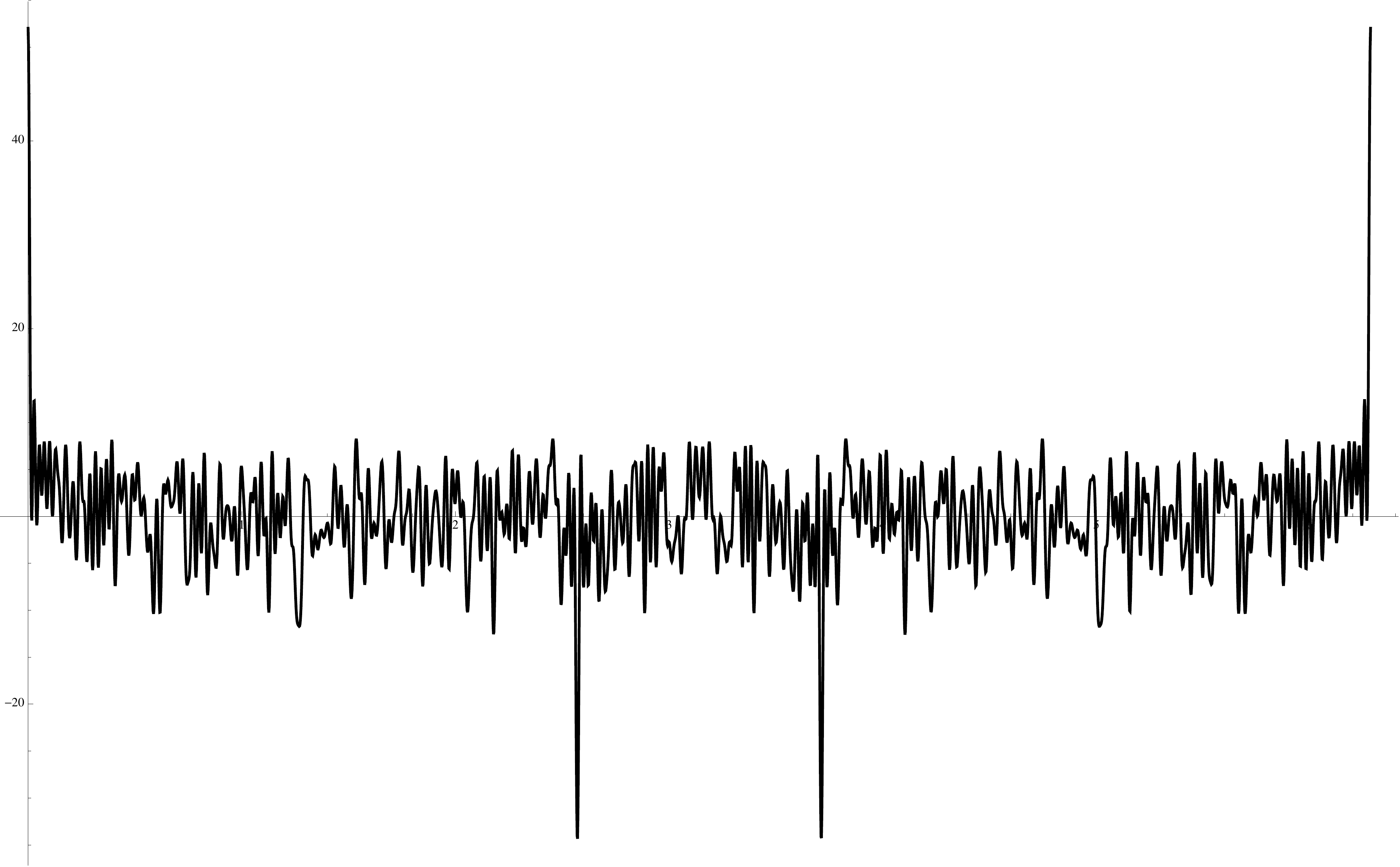}
\caption{The function $f_{50}$.}
\end{center}
\end{minipage}%
\begin{minipage}[t]{0.49\columnwidth}%
\begin{center}
\includegraphics[width = 5cm]{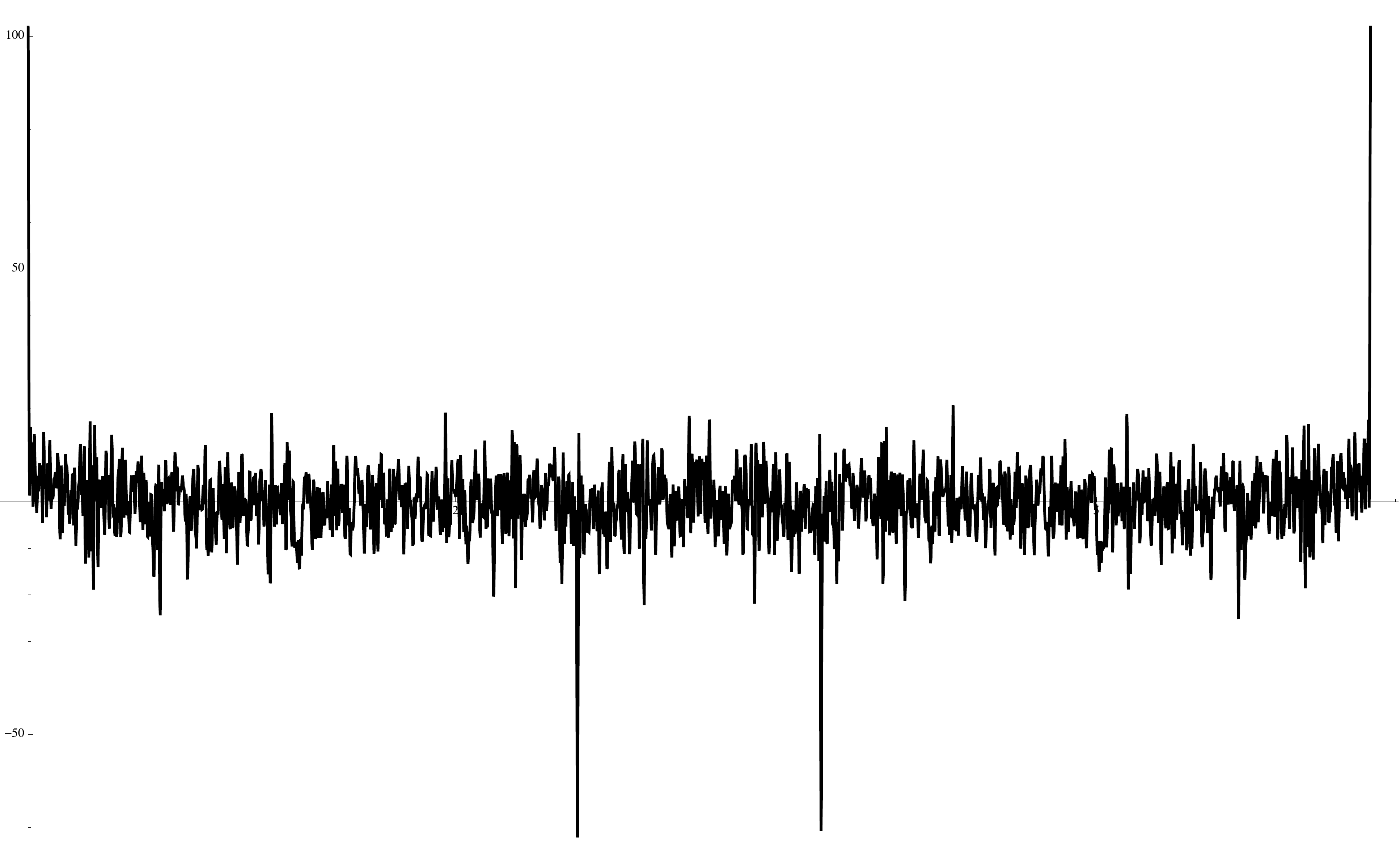}
\caption{The function $f_{100}$.}
\end{center}
\end{minipage}
\end{figure}

\subsection{The hidden signal.} An explicit computation with $N = 10^7$ pinpoints the location of the signal at
$$ \alpha \sim 2.5714474995\dots$$
and, by symmetry, at $2\pi - \alpha$. This signal acts as a hidden shift in frequency space: we can remove that shift in frequency by considering
$$S_N= \left\{ \alpha a_n - 2\pi \left\lfloor \frac{\alpha a_n}{2\pi} \right\rfloor: 1 \leq n \leq N \right\} \qquad \mbox{instead}.$$
\begin{figure}[h!]
\begin{center} 
\includegraphics[width = 10cm]{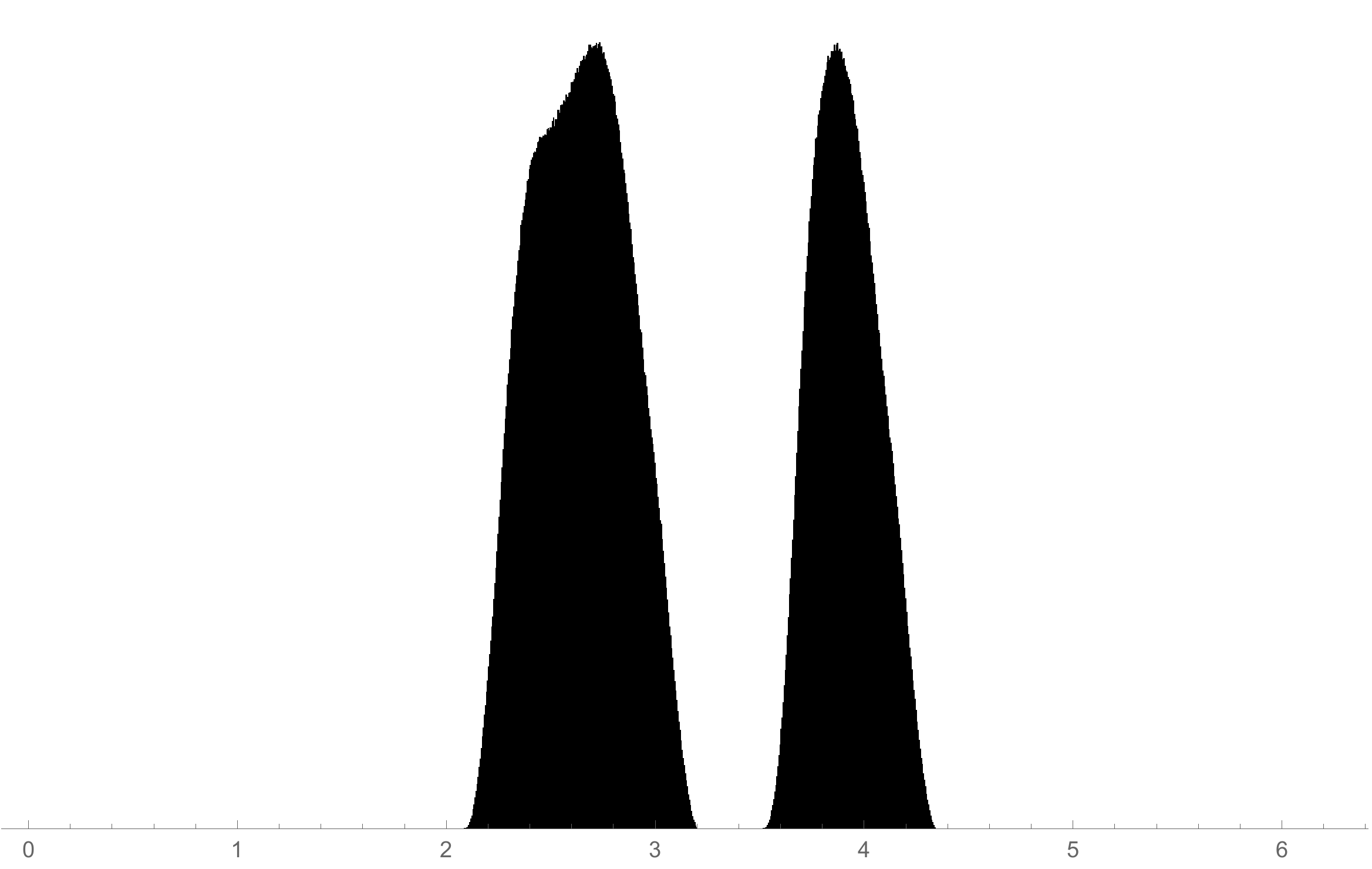}
\caption{Distribution of $S_N$ on $[0,2\pi]$ for $N=10^7$.}
\end{center}
\end{figure}

A priori it might be reasonable to suspect that the Ulam sequence is `pseudo-random' (in the same way as Cram\'{e}r's model suggests `pseudo-randomness' of the primes). There are some local obstructions (for example both primes and the Ulam sequence do not contain a triple $(n, n+2, n+4)$ for $n > 3$) but it would a priori be conceivable that many statistical quantities could be predicted by a random model. Our discovery clearly shows this to be false: the Ulam sequence has some \textit{extremely} rigid underlying structure. We are hopeful that the rigidity of the phenomenon will provide a first way to get a more
rigorous understanding of the Ulam sequence; however, we also believe that the phenomenon might be related to some interesting mechanism in additive combinatorics and of independent interest. 

\subsection{$\alpha$ and the power spectrum.} Gibbs \cite{gibbs}  has used this phenomenon as a basis for an algorithm that allows for a faster
computation of Ulam numbers. Jud McCranie informed me that he and Gibbs have independently computed a large number
of Ulam numbers (more than 2.4 billion) and verified the existence of the phenomenon up to that order. McCranie gives the bounds
$$ 2.57144749846 < \alpha  < 2.57144749850.$$
We do not have any natural conjecture for a closed-form expression of $\alpha$. The value $\alpha$ is not unique in the sense that there are other values $x$ such that
$$ \sum_{n=1}^{N}{\cos{(x a_n)}} \sim c_{x}N \qquad \mbox{seems to exhibit linear growth,}$$
however, these other values can be regarded as the effect of an underlying symmetry.
Sinan G\"unt\"urk (personal communication) observed that these other values seem to be explicitly given by $\left\{k\alpha~\mbox{mod}~2\pi: k \in \mathbb{Z}\right\}$. If there is weak convergence of the empirical distribution, i.e. if
$$  \frac{1}{N} \sum_{n=1}^{N}{\delta_{\alpha a_n}} \rightharpoonup \mu_{\alpha}~\mbox{on}~\mathbb{T}\mbox{, then} \qquad \lim_{N \rightarrow \infty}{\frac{1}{N}\sum_{n=1}^{N}{\cos{(\alpha a_n)}}} = \int_{\mathbb{T}}{\cos{(x)} d\mu_{\alpha}}.$$
Computationally, our strategy
consists of looking for nonzero values of the integral (whose value is 0 in the generic case, where $\mu_{\alpha}$ is a uniform distribution). In this framework, the frequency $\alpha \sim 2.571\dots$ is the one for which the phenomenon is most pronounced, however, there are other peaks. It is easy to see that if we have convergence to an absolutely continuous measure
$$  \frac{1}{N} \sum_{n=1}^{N}{\delta_{\alpha a_n}}  \rightharpoonup f(x)dx,~\mbox{then for all}~\ell \in \mathbb{Z} \qquad \frac{1}{N}\sum_{n=1}^{N}{\delta_{\alpha \ell a_n}}  \rightharpoonup \left( \frac{1}{\ell} \sum_{k=0}^{\ell - 1}{ f\left(\frac{x+2 \pi k}{\ell}\right)}\right)dx$$
and if $\mu = f(x)dx$ is unevenly distributed in the sense that
$$ 0 \neq  \int_{\mathbb{T}}{\cos{(x)} d\mu} =  \int_{\mathbb{T}}{\cos{(x)} f(x) dx},$$ 
then the same relation can be expected for the new density function. There are functions
$$ f \neq \frac{1}{2\pi} \quad \mbox{for which nonetheless} \quad \frac{1}{\ell}\sum_{k=0}^{\ell - 1}{ f\left(\frac{x+2 \pi k}{\ell}\right)} \equiv \frac{1}{2\pi},$$
however, the set of such functions is very small and not stable under perturbations. \textit{Generically}, if $f(x)dx$ is not
the uniform distribution, then neither is the new distribution. In the case of the Ulam sequence $(a_n)$ it is thus to be expected that the signals at frequency $\ell\alpha$ mod $2\pi$ generated out of the initial measure (Fig.3) located at frequency $\alpha$ are never uniformly distributed. Numerically, this holds and
$$ \sum_{n=1}^{N}{\cos{(\ell \alpha a_n )}} \sim c_{\ell} N.$$
\vspace{-5pt}
\begin{table}[h!] 
\begin{tabular}{ | l |c | c | c | c | c | c | c | c | c | c }
\hline	
$\ell$ & $0$ & $1$ & $2$ & $3$ & $4$ & $5$  & $6$ & $7$ & $8$  	\\
\hline 
$c_{\ell}$ & 1 & -0.794 & 0.288 & 0.253 & -0.578 & 0.580 & -0.344 & 0.057 & 0.118  \\	
\hline
\end{tabular} 
\\[5pt]
 \captionsetup{width=.9\linewidth}
\caption{Empirical approximations of $c_{\ell}$ (and $c_{-\ell} = c_{\ell}$).}
\end{table}

Numerical computation suggests that indeed \textit{all} peaks in the spectrum are described by the set $\left\{k\alpha~\mbox{mod}~2\pi: k \in \mathbb{Z}\right\}$ which suggests that there is indeed only one hidden signal (shown in
Figure 3) at frequency $\alpha \sim 2.571\dots$ while all other peaks are explained by the symmetry described above. 

\begin{figure}[h!]
\begin{minipage}[t]{0.49\columnwidth}%
\begin{center}
\includegraphics[width = 7cm]{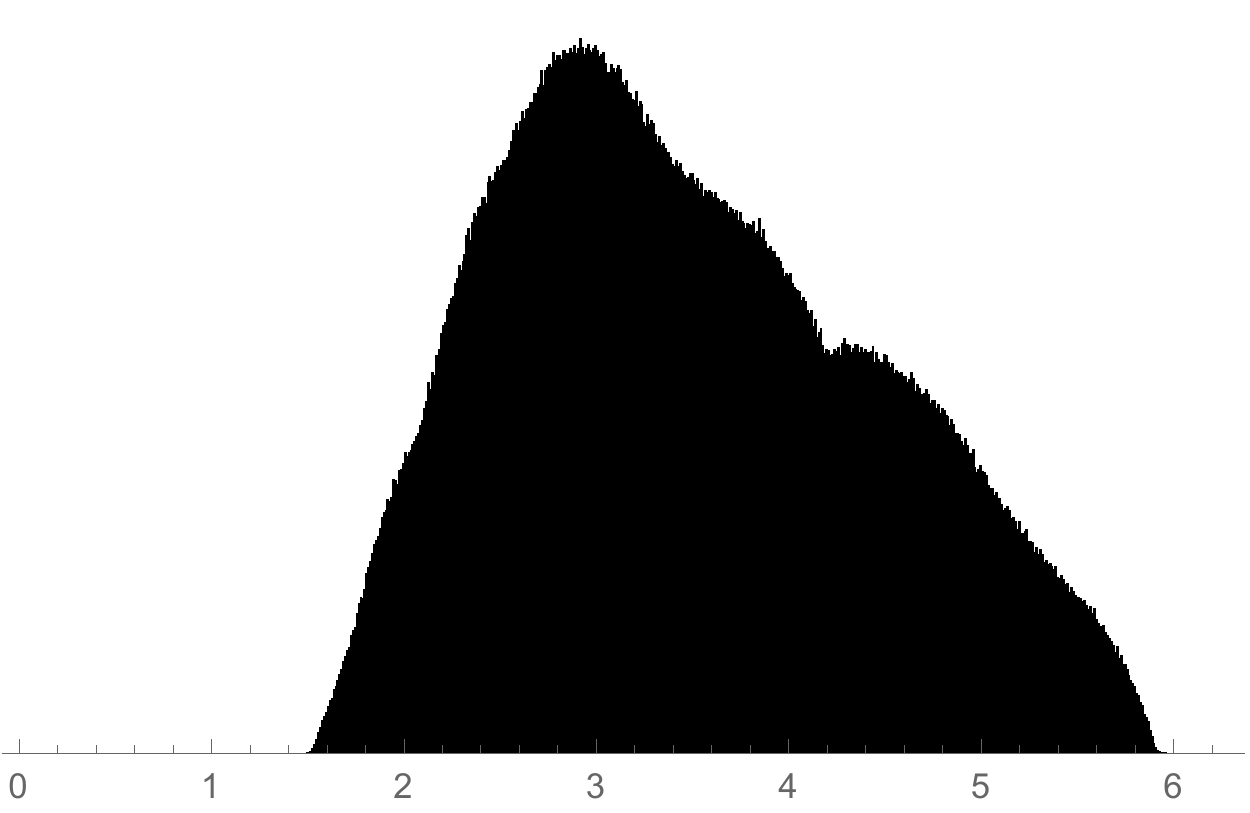}
\captionsetup{width=0.9\textwidth}
\caption{Distribution at frequency $4\alpha$.}

\end{center}
\end{minipage}%
\begin{minipage}[t]{0.49\columnwidth}%
\begin{center}
\includegraphics[width = 7cm]{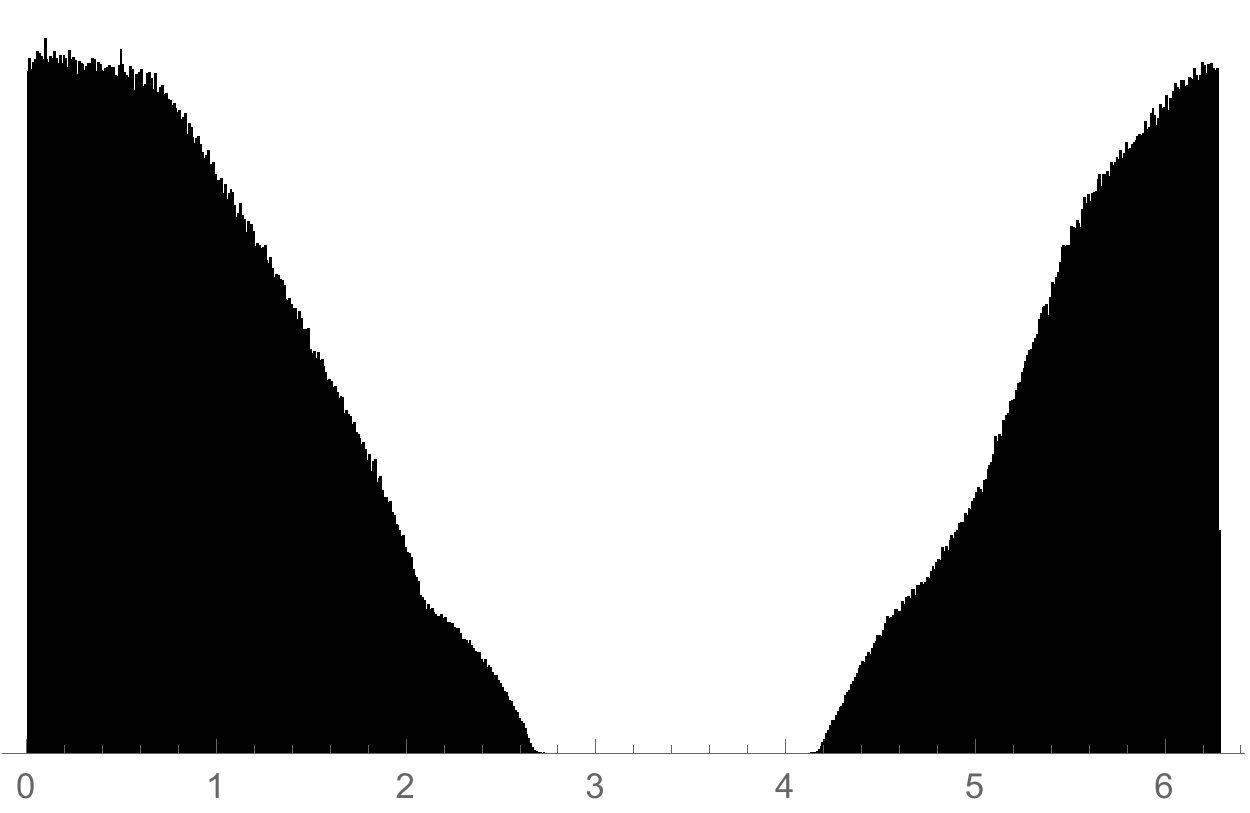}
\captionsetup{width=0.9\textwidth}
\caption{Distribution at frequency $5\alpha$.}
\end{center}
\end{minipage}
\end{figure}

 $c_4$ and $c_5$ are especially large compared to other values (though not as big as $c_1 \sim 2.571\dots$). 
This seems to be without deeper meaning: the distributions arising for those values are mainly supported in regions
where the cosine is negative and positive, respectively. It is easy to see that under suitable assumptions on the initial measure we get that $c_{\ell} \rightarrow 0$ as
$\ell \rightarrow \infty$.

\subsection{Other initial values}
There has been quite some work on the behavior of Ulam-type sequences using the same construction rule
with other initial values. It has first been observed by Queneau \cite{q} that
the initial values $(a_1, a_2) = (2,5)$ give rise to a sequence where $a_{n+1} - a_{n}$ is eventually periodic.
Finch \cite{finch1, finch2, finch3} proved that this is the case whenever only finitely many even numbers appear in the
sequence and conjectured that this is the case for $(a_1, a_2) = (2,n)$ whenever $n \geq 5$ is odd.
Finch's conjecture was proved by Schmerl \& Spiegel \cite{schm}.  Subsequently, Cassaigne \& Finch \cite{cas}
proved that all sequences starting from $(a_1, a_2) = (4,n)$ with $n \equiv 1$ (mod 4) contain precisely three even integers and are thus eventually periodic. The
sequence $(a_1, a_2) = (2,3)$ as well as all sequences $(a_1, a_2) = (1,n)$ with $n \in \mathbb{N}$
do not exhibit such behavior and are being described 'erratic' in the above literature. 

\begin{figure}[h!]
\begin{minipage}[t]{0.49\columnwidth}%
\begin{center}
\includegraphics[width = 7cm]{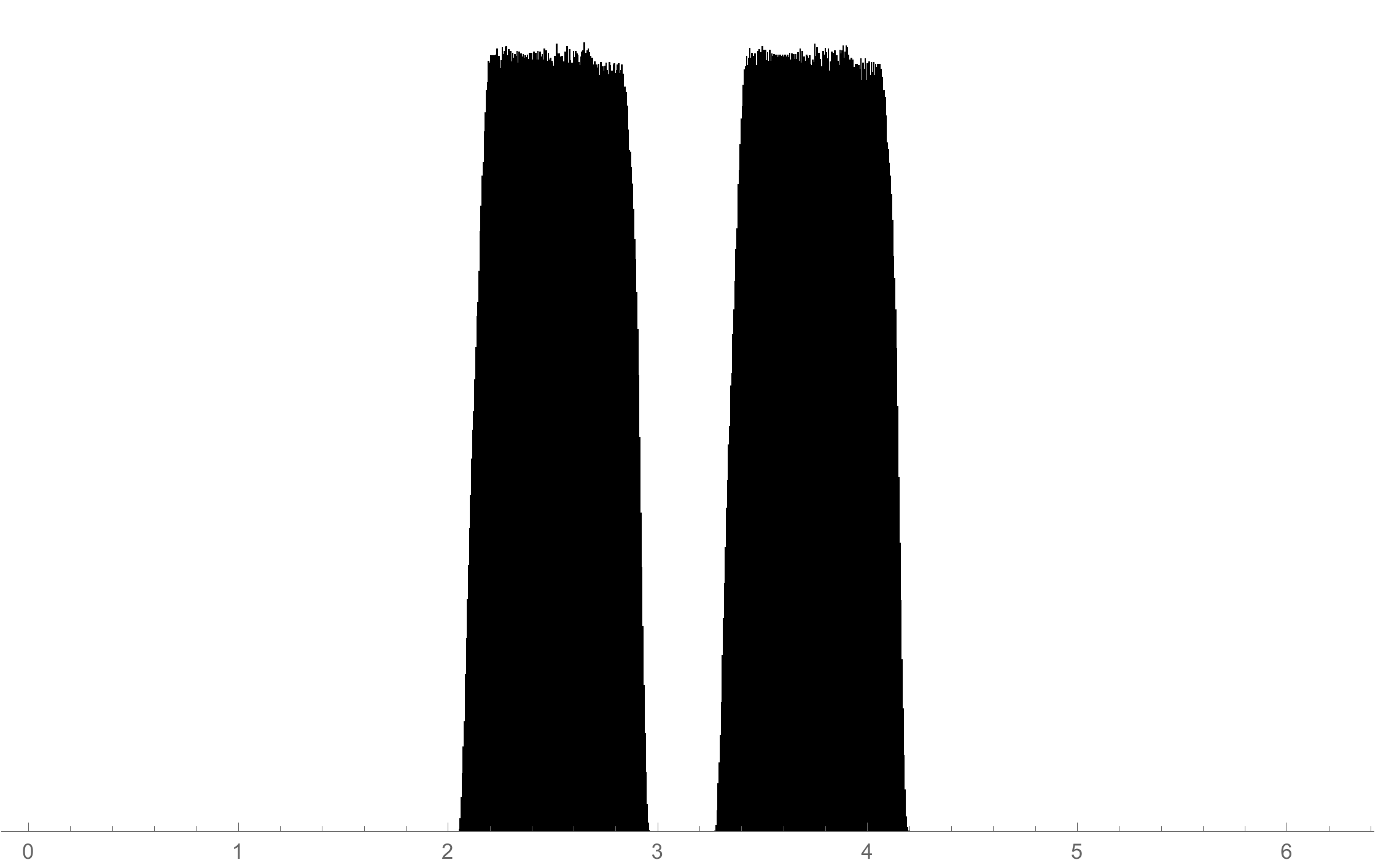}
\captionsetup{width=0.9\textwidth}
\caption{Distribution for initial values $(a_1, a_2) = (1,3)$ for $N = 2.5 \cdot 10^6$.}

\end{center}
\end{minipage}%
\begin{minipage}[t]{0.49\columnwidth}%
\begin{center}
\includegraphics[width = 7cm]{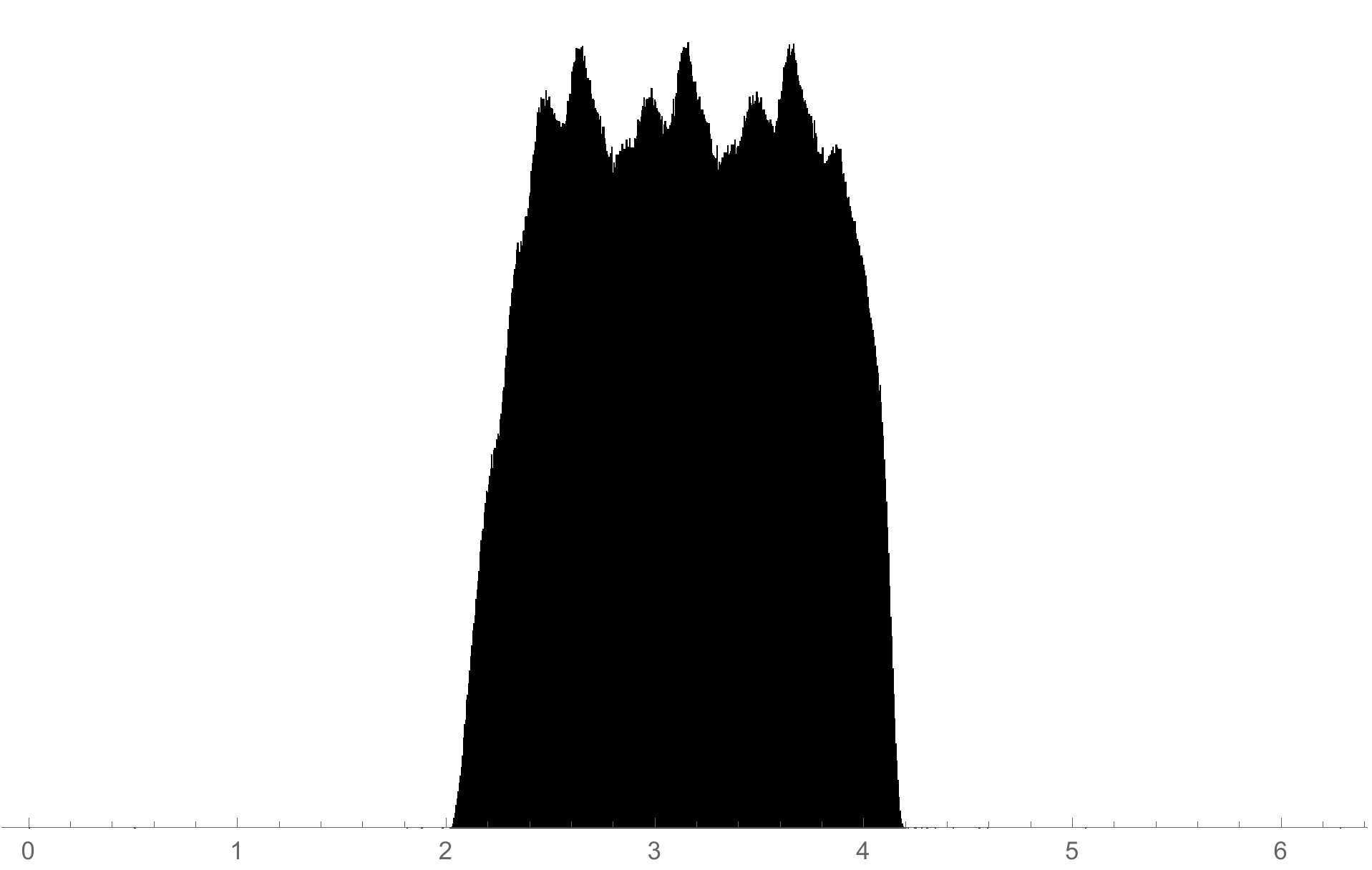}
\captionsetup{width=0.9\textwidth}
\caption{Distribution for the initial values $(a_1, a_2) = (1,4)$ for $N = 3.9 \cdot 10^6$.}
\end{center}
\end{minipage}
\end{figure}

The underlying mechanism appears to be independent of the initial values (as long as the sequence is 'chaotic' in the sense of not having periodic consecutive differences). Using $N=2.7 \cdot 10^6$ and $N = 3.9 \cdot 10^6$ elements, respectively, we found that the frequencies for the erratic sequences with $(a_1, a_2) = (1,3)$ and $(a_1, a_2) = (1,4)$ are given by
$$\alpha_{(1,3)} \sim 2.83349751\dots  \qquad \mbox{and} \qquad \alpha_{(1,4)} \sim 0.506013502\dots$$
and removing that hidden frequency gives rise to the two distributions Fig. 5 and Fig. 6. 
Another example is given by initial conditions $(2,3)$ with $N = 5.7 \cdot 10^6$ elements,
$$\alpha_{(2,3)} \sim 1.1650128748\dots$$
\vspace{-10pt}
\begin{figure}[h!]
\begin{center}
\includegraphics[width = 9cm]{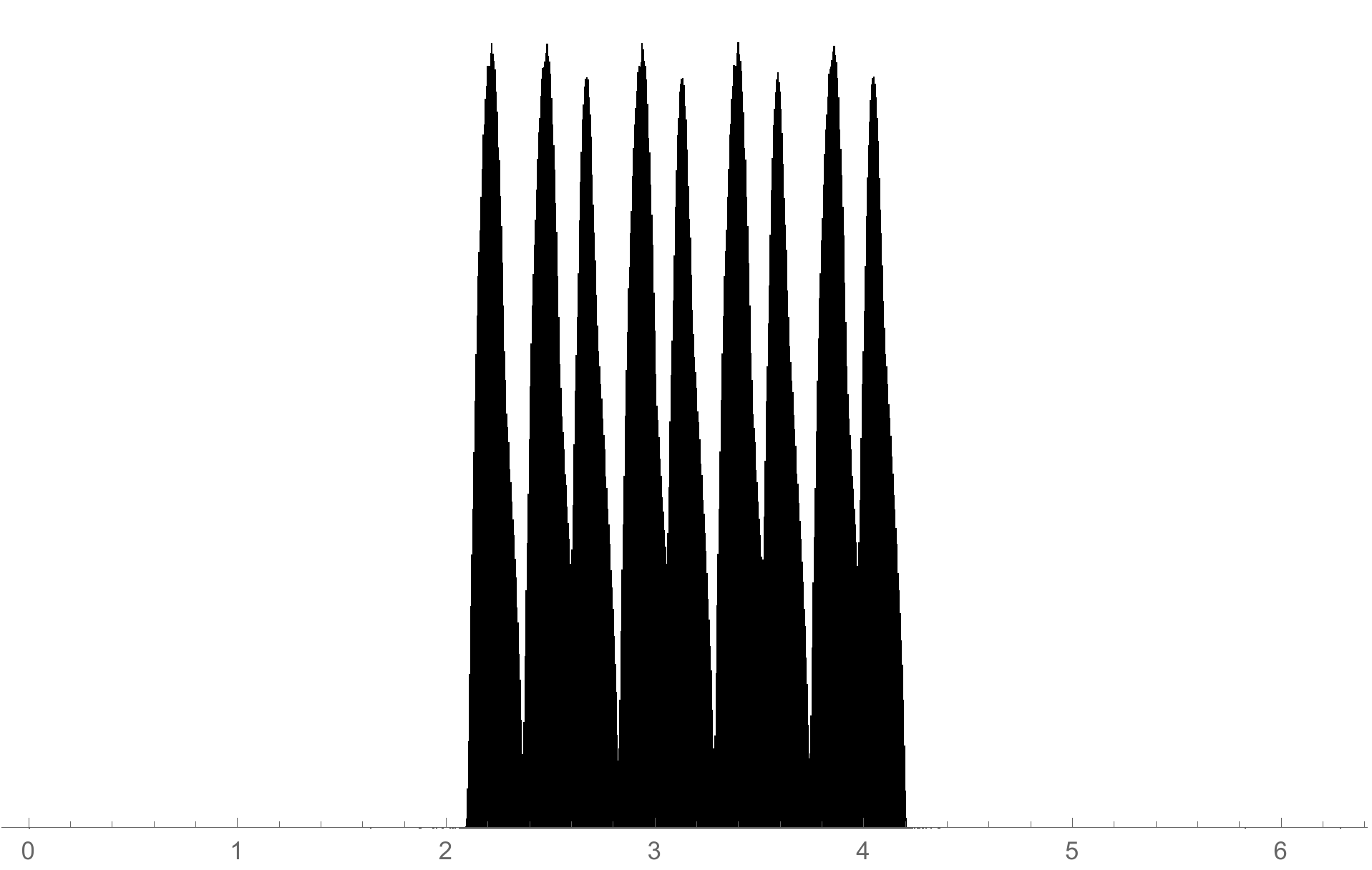}
\caption{Distribution for the initial values $(2,3)$ for $N=5.7 \cdot 10^6$.}
\end{center}
\end{figure}
\vspace{-10pt}
\subsection{Additional remarks} 
 A classical theorem of Hermann Weyl \cite{weyl} states that if $(a_n)_{n=1}^{\infty}$ is a sequence of distinct positive integers, then
the set of $\alpha \in \mathbb{R}$ for which the sequence 
$$(\alpha a_n~\mbox{mod}~2\pi)_{n=1}^{\infty}~\mbox{is \textit{not} uniformly distributed}$$
has measure 0. Conversely, given any absolutely continuous measure $\mu$ on $\mathbb{T}$ and a fixed $\alpha$ such that $\alpha/(2\pi)$ is irrational, it is not difficult to construct a sequence of integers $a_n$ such that 
 $(\alpha a_n~\mbox{mod}~2\pi)_{n=1}^{\infty}$ is distributed according to $\mu$ by using the uniform distribution of $(k\alpha)_{k \in \mathbb{N}}$. Naturally, these sequences are non-deterministic and quite artificial. 
\begin{quote}
\textbf{Question.} Are there any other sequences $(a_n)_{n \in \mathbb{N}}$ of integers appearing 'naturally' with the property that there exists
a real $\alpha > 0$ such that
$$(\alpha a_n~\mbox{mod}~2\pi)_{n=1}^{\infty}~\mbox{has an absolutely continuous \textit{non-uniform} distribution?}$$
\end{quote}
We expect this property to be exceedingly rare.

\section{Other examples of curious behavior}
The author, after discovering the Ulam sequence phenomenon, tried a very large number of different sequences arising in combinatorics and number theory in the hope of uncovering a second example. We have not found any other sequence with that property. Generally, a natural case distinction was observed:
\begin{enumerate}
\item $\cos{(a_1 x)} + \cos{(a_2 x)} + \dots + \cos{(a_N x)}$ appears random
\item $\cos{(a_1 x)} + \cos{(a_2 x)} + \dots + \cos{(a_N x)}$ has peaks at $2\pi/k$ for some $k \in \mathbb{N}$
\item $\cos{(a_1 x)} + \cos{(a_2 x)} + \dots + \cos{(a_N x)}$ exhibits other behavior.
\end{enumerate}
The first type of behavior is generic and observed for most sequences.
The second type of behavior merely indicates irregularities in the distribution of the sequence modulo $k$ and is
easily observed in many sequences. We present one known example (based on work of Reznick \cite{rez}.
The third type seems more complicated and it is not clear
what type of phenomena one could encounter. In trying a very large number of examples, we only found the following three examples and
they appear to be very different from each other:
\begin{itemize}
\item a fixed, stationary distribution (Ulam sequence)
\item a fixed distribution on top of the uniform distribution, which weakens as more elements
are added (the zeroes of the Riemann $\zeta$-function on the critical line, explained
by Landau \cite{landau}, Ford \& Zaharescu \cite{ford})
\item different intensities of randomness (Lagarias' modification of 'primes of measurement')
\end{itemize}

\subsection{Example 1: Stern diatomic sequence (Dijkstra, Reznick).} We start with an interesting example of a sequence having irregulaties in residue classes. The Stern diatomic sequence is defined by $a_0 = 0$, $a_1 = 1$ and
$$ a_{n} = \begin{cases} a_{n/2} \qquad &\mbox{if}~n~\mbox{is even} \\
a_{\frac{n-1}{2}} + a_{\frac{n+1}{2}} \qquad &\mbox{if}~n~\mbox{is odd.} 
\end{cases}$$
The sequence has a surprising number of combinatorial properties. Plotting the associated cosine sum easily identifies peaks (of decreasing
height) at $2\pi/3, 2\pi/5, 2\pi/7, \dots$ -- it is easy to find the origin of these peaks: elements of the Stern sequence are
equally likely to be $\equiv i $ (mod $p$) as long as $i \neq 0$ but they are slightly less likely to be divisible by $p$. The simplest possible
case ($2 \big| a_n$ if and only if $3 |n$) was observed by Dijsktra \cite{dij} in 1976 (with very similar arguments being already contained
in the original paper of Stern \cite{stern}). The fact that the phenomenon persists mod $p \neq 2$ seems to have first been
discovered and proven by Reznick \cite{rez} in 2008. 
This example is representative: typically, if there any peaks, they are caused by an irregular distribution mod $p$.

\begin{figure}[h!]
\begin{minipage}[t]{0.49\columnwidth}%
\begin{center}
\includegraphics[width = 6cm]{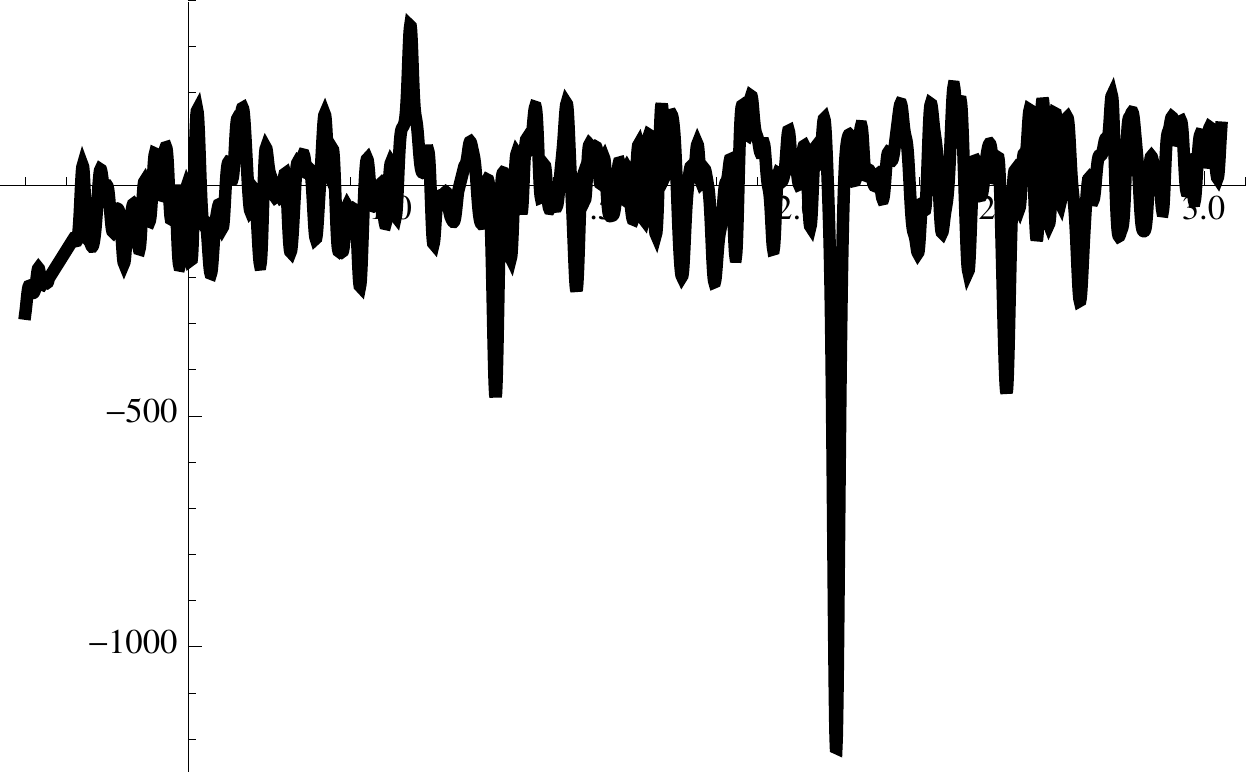}
\captionsetup{width=0.9\textwidth}
\caption{Distribution for $N = 10^4$.}

\end{center}
\end{minipage}%
\begin{minipage}[t]{0.49\columnwidth}%
\begin{center}
\includegraphics[width = 6cm]{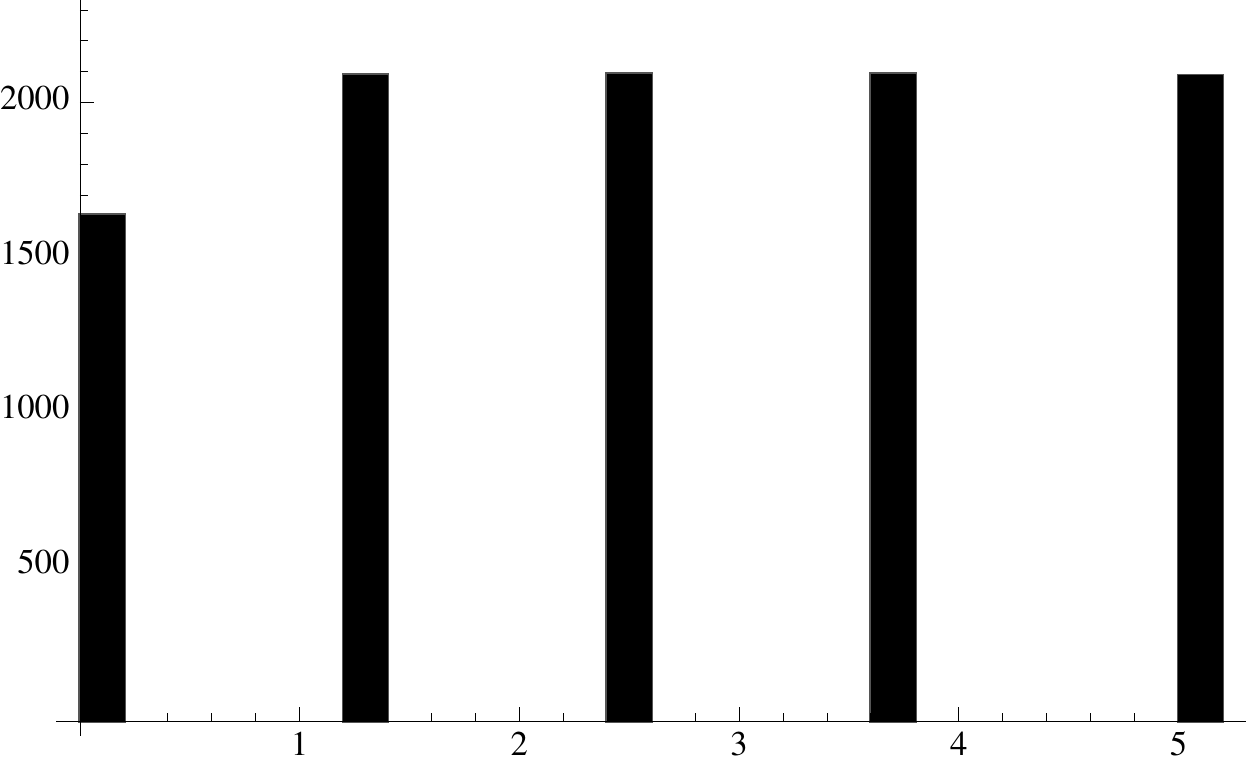}
\captionsetup{width=0.9\textwidth}
\caption{$\alpha a_n ~\mbox{mod}~2\pi$ with $\alpha = 2\pi/5$.}
\end{center}
\end{minipage}
\end{figure}

\subsection{Example 2: Zeroes of the Riemann $\zeta-$function (Landau, Ford-Zaharescu).} This section does not contain any new results but describes known result in our framework. We now discuss an interesting example of a sequence that exhibits a
fixed distribution that weakens as more elements are added:
the sequence $t_n$ of imaginary parts of zeroes $\zeta(1/2 + i t_n) = 0$ of the Riemann zeta function on the critical line. The $\zeta-$function
is defined via
$$ \zeta(s) = \sum_{n \geq 1}{\frac{1}{n^s}} \qquad \mbox{and we are interested in} \qquad \zeta\left(\frac{1}{2} + i t \right), ~t \in \mathbb{R}.$$
The first few zeroes $(t_n)$ are roughly located at
$$ \sim 14.13,~ 21.02, ~25.01, ~ 30.42, ~ 32.93, ~ 37.58, \dots$$
It is well-known that
$$ \#\left\{t \in [0,T]: \zeta\left(\frac{1}{2} + i t \right) = 0 \right\} \sim \frac{T \log{T}}{2 \pi} \qquad \mbox{as} \quad T \rightarrow \infty.$$
This means that the density of zeroes $t_1 \leq t_2 \leq t_3 \leq \dots$ on the critical line is increasing and
their consecutive differences are shrinking like $1/\log{T}$. The distribution of consecutive differences of properly renormalized roots has been
intensively studied ever since Montgomery \cite{mont} discovered a connection to distributions arising in the Gaussian Unitary Ensemble (GUE).
This is of great relevance in light of the Hilbert-Polya question whether the roots correspond to the eigenvalues of a self-adjoint linear operator. 
\begin{figure}[h!]
\begin{minipage}[t]{0.49\columnwidth}%
\begin{center}
\includegraphics[width = 7cm]{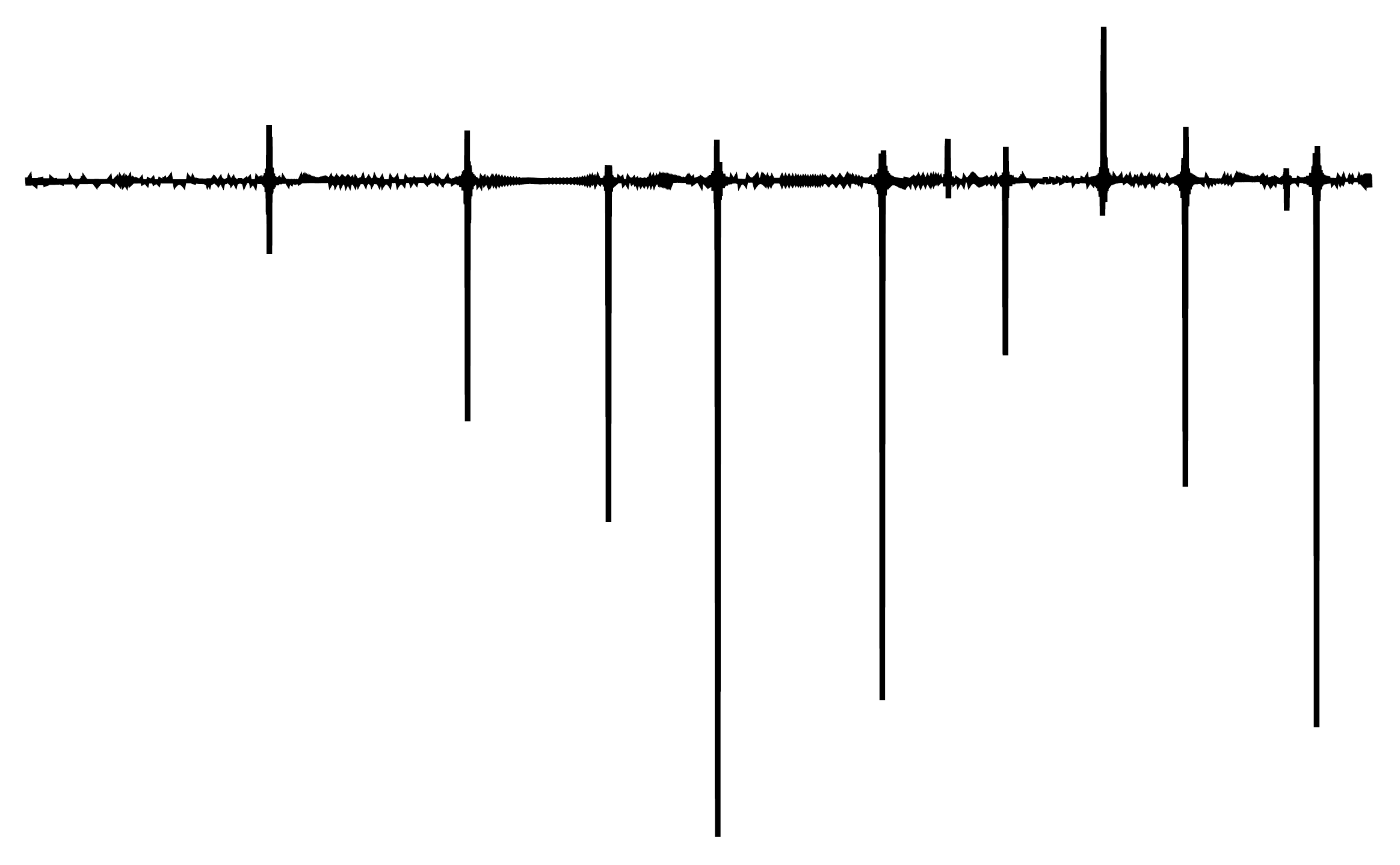}
\captionsetup{width=0.9\textwidth}
\caption{Value of the cosine sum on the interval $0.1 \leq x \leq \pi - 0.1$ using $10^5$ roots.}

\end{center}
\end{minipage}%
\begin{minipage}[t]{0.49\columnwidth}%
\begin{center}
\includegraphics[width = 7cm]{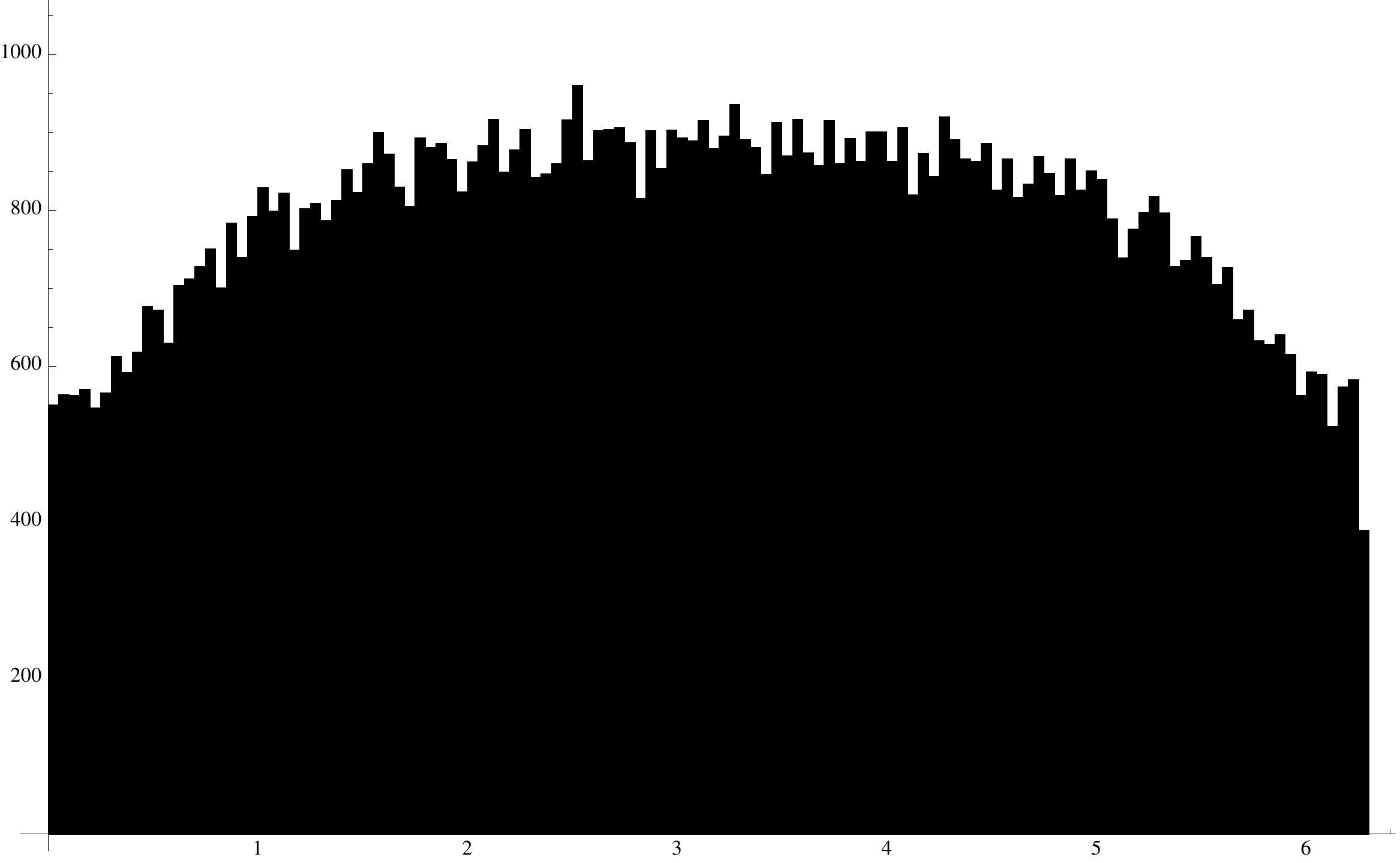}
\captionsetup{width=0.9\textwidth}
\caption{The distribution of $(\log{5}) t_n ~\mbox{mod}~2\pi$ for the first $10^5$ roots.}
\end{center}
\end{minipage}
\end{figure}
Computing the cosine sum using the zeroes as frequencies results in a most interesting picture (Fig. 11); the spikes are centered at $\log{p^m}$ for $p$ prime and $m \in \mathbb{N}$. 
This is a consequence of a result of Landau \cite{landau}, who proved (\textit{without} assuming the Riemann hypothesis)
$$ \sum_{n=1}^{N}{e^{i t_n x}} = \begin{cases} - \frac{t_N}{2\pi} \frac{\log{p}}{\sqrt{\log{p^m}}}  + \mathcal{O}\left(  \frac{\log{N}}{\sqrt{\log{p^m}}} \right)  \qquad &\mbox{if}~x = \log{p^m} \\
\mathcal{O}\left(  \frac{\log{N}}{\sqrt{\log{p^m}}} \right)  \qquad &\mbox{otherwise.}
\end{cases}$$ 
Note that the sum ranges over $N$ elements and the leading order term for $x = \log{p^m}$ is given by $t_N \sim N/\log{N}$. In order for this to happen at least $N/\log{N}$ out of the first $N$ zeroes of the zeta function have to align in a nontrivial way when considered as $\log{(p^m)}t_n ~\mbox{mod}~2\pi$ (see Table 1).
Since $N/\log{N} \ll N$ that alignment weakens as more zeroes are being added. The distribution at scale 
$N/\log{N}$ was found by Ford \& Zaharescu \cite{ford}.\\

There is a simple heuristic for this clustering using Euler's product formula
$$ \zeta(s) = \sum_{n=1}^{\infty}{\frac{1}{n^s}} = \prod_{p}{\frac{1}{1-p^{-s}}}.$$
$\zeta(s) = 0$ means that the infinite product converges to 0 (no factors vanishes), which requires
the prime numbers to be suitably aligned. The relationship with the observed alignment is due to 
$$ \left| \frac{1}{1-p^{-\left(\frac{1}{2} + it_n\right)}} \right| \leq 1 \qquad \mbox{if} \qquad \frac{\pi}{2} \leq \left[ (\log{p}) t_n~\mbox{mod}~2\pi \right] \leq \frac{3 \pi}{2}.$$
\vspace{-5pt}
\begin{table}[h!] 
\begin{tabular}{ | l |c | c | c | c | c | c | c | c | c | c }
\hline	
$\alpha$ & $\log{2}$ & $\log{3}$ & $\log{5}$ & $\log{7}$ & $\log{8}$ & $\log{9}$  & $\log{11}$ & $e$ & $\sqrt{5}$ 	\\
\hline 
$\#$ points & 53258 &54392 & 55123 & 55336 & 51883 & 52572 & 55398 & 49992 & 50086  \\	
\hline
\end{tabular} 
\\[6pt]
\caption{Number of elements in $\left\{\alpha t_i ~\mbox{mod}~2\pi: 1 \leq i \leq 10^5 \right\}$ that lie in $[\pi/2, 3\pi/2]$.}
\end{table}
\vspace{-20pt}

\subsection{Example 3: Lagarias' variant of MacMahon's primes of measurement.} MacMahon's 'primes of measurement', sometimes called segmented numbers, are given by the sequence
$$ 1,2,4,5,8,10,14,15,16,21,22,25,26,28, \dots$$
generated by excluding all sums of two or more consecutive earlier elements of the sequence (see Guy's \textit{Unsolved Problems in Number Theory} \cite[Section E30]{guy}). George
Andrews conjectures that
$$ a_n \sim \frac{n \log{n}}{\log \log n}.$$
The sequence is very different from the Ulam sequence but certainly similar in spirit. The cosine sum with these terms does appear to be random, however, we also remark that the publicly available data set \cite{oeis3}
we used only contains about 7000 terms and larger compilations may provide more insight. Lagarias \cite{lagarias} proposes excluding only sums of two or three consecutive earlier members
(which certainly increases the similarity to the definition of the Ulam sequence), the arising sequence being
$$ 1,2,4,5,8,10,12,14,15,16,19,20,21,24,25,27,28, \dots$$
We used a publicly available list of the first 10.000 elements \cite{oeis4}. First of all we note a peak at $\pi$ corresponding to an uneven distribution in $\mathbb{Z}_{2}$: indeed, it seems that elements of the sequence are more likely to be even than odd (with $\sim 54.86\%$ of elements being even).

\begin{figure}[h!]
\begin{minipage}[t]{0.49\columnwidth}%
\begin{center}
\includegraphics[width = 7cm]{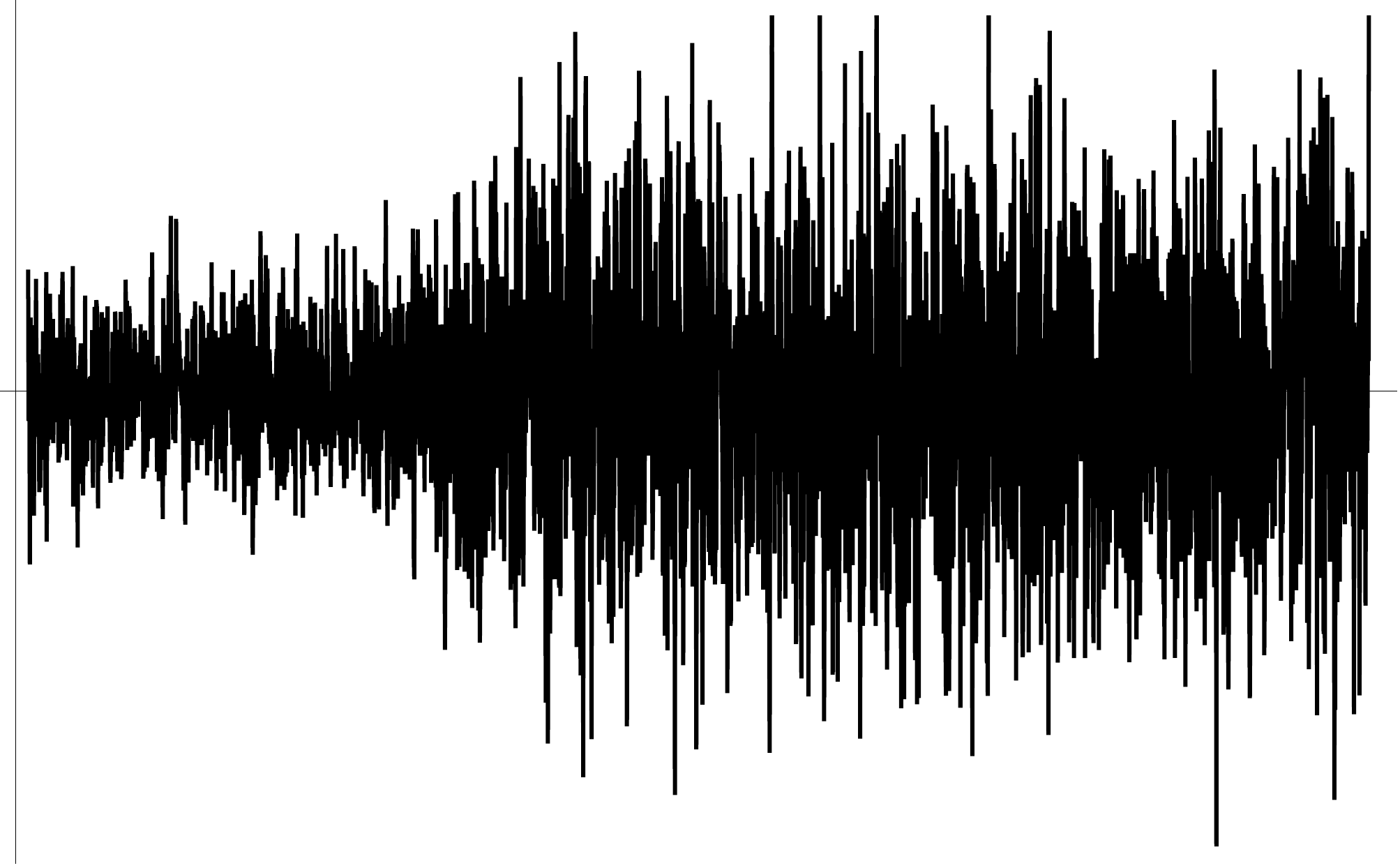}
\end{center}
\end{minipage}%
\begin{minipage}[t]{0.49\columnwidth}%
\begin{center}
\includegraphics[width = 7cm]{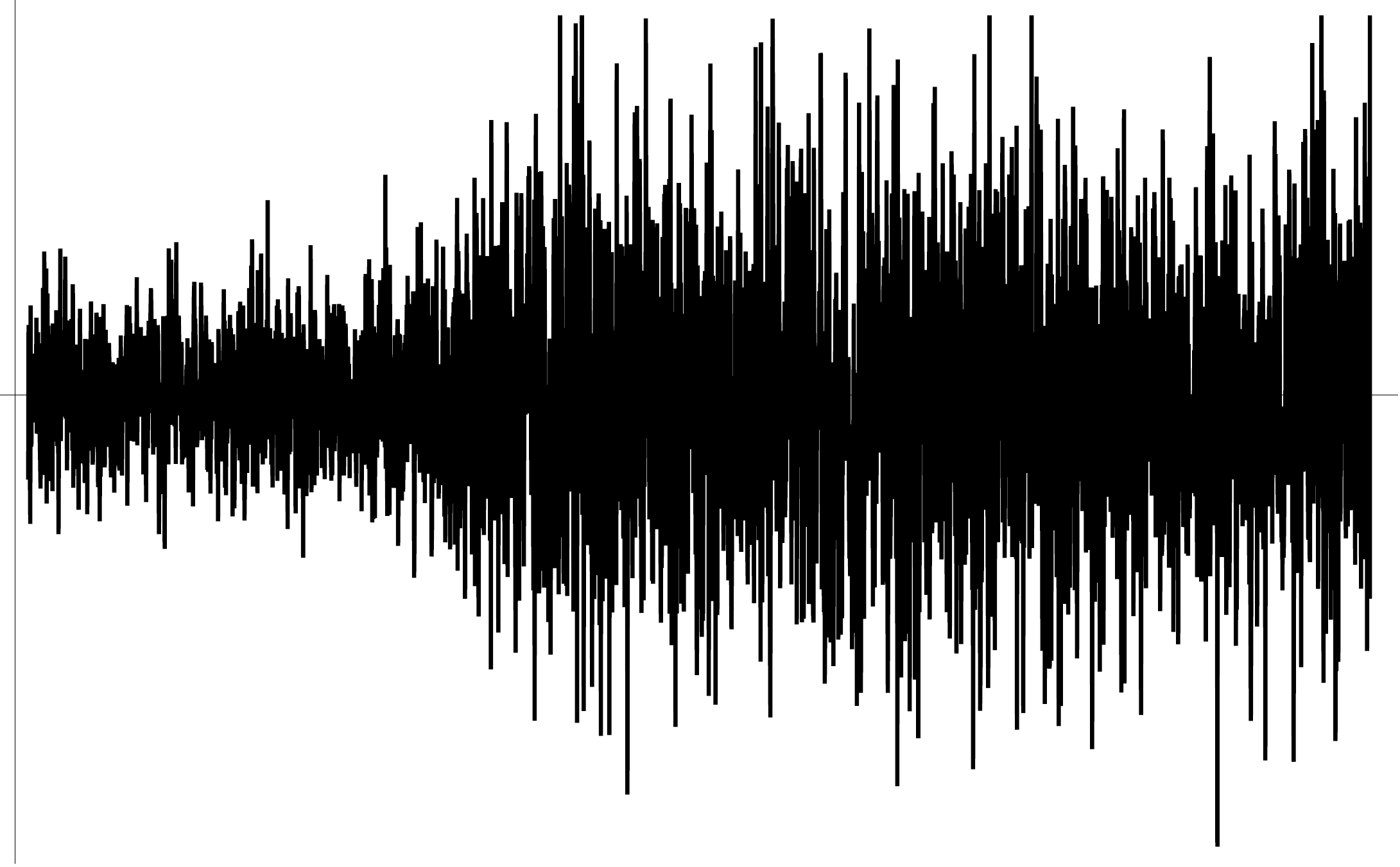}
\end{center}
\end{minipage}
\captionsetup{width=0.9\textwidth}
\caption{Cosine sum on $[0.03, \pi-0.03]$ for the first $5.000$ and $10.000$ elements.}
\end{figure}

The behavior away from the peak is more interesting: it seems as if 
$$ \left|  \sum_{n=1}^{N}{\cos{(a_n x)}} \right| \leq c \sqrt{N} \qquad \mbox{for $x$ away from 0}$$
but that the sum is experiencing smaller levels of fluctuations for small $x$ and larger levels for large $x$. There seems to be a transition occuring somewhere
around $x \sim 1$ . The absolute value of the sum seems to be indeed of order $\sim \sqrt{N}$ and, as a consequence, plotting the distribution of $\left\{x a_n~\mbox{mod}~2\pi\right\}$
on the torus shows essentially a uniform distribution.

\medskip
\medskip

\textbf{Acknowledgement.} The result for the Ulam sequence was originally discovered using only the first $10.000$ numbers for each set of initial conditions; the precision of the results 
presented here would not have been possible without the data sets compiled and generously provided by Daniel Strottman for all initial conditions that were discussed. 
Sinan G\"unt\"urk observed the arithmetic regularity of the location of the peaks which
gave rise to Section 2.3. Bruce Reznick was very helpful in explaining the early history of the Stern sequence.
Jud McCranie informed the author of additional data and performed additional tests on them providing a better estimate for $\alpha_{(1,2)}$. The computations 
involving the $\zeta-$function were carried out using Oldyzko's list \cite{od2} of the first 100.000 roots. The author 
is indebted to Steven Finch for extensive discussions and his encouragement.

\end{document}